\documentclass[10pt,rotate]{amsart}

\usepackage[all]{xy}

\usepackage{latexsym}

\usepackage{amssymb}
\usepackage{amsfonts}
\usepackage{amscd}
\usepackage{amsmath,amsthm}

\newtheorem{lemma}{Lemma}[section]

\newtheorem{lem}[lemma]{Lemma}
\newtheorem{prop}[lemma]{Proposition}
\newtheorem{thm}[lemma]{Theorem}
\newtheorem{cor}[lemma]{Corollary}

{

}

\theoremstyle{definition}

\newtheorem{example}[lemma]{Example}

\theoremstyle{remark}

\numberwithin{equation}{section}

\newenvironment{pf}{\noindent{\bf Proof.}}{\hfill $\square$\medskip}

\def\PP{{\mathbb P}}

\def\ZZ{{\mathbb Z}}

\def\0ol{{\bar 0}}
\def\1ol{{\bar 1}}
\def\2ol{{\bar 2}}
\def\ol2{{\bar 2}}
\def\3ol{{\bar 3}}
\def\4ol{{\bar 4}}
\def\5ol{{\bar 5}}
\def\6ol{{\bar 6}}
\def\7ol{{\bar 7}}
\def\8ol{{\bar 8}}
\def\9ol{{\bar 9}}

\def\bold0{{\bf 0}}
\def\bold1{{\bf 1}}
\def\bold2{{\bf 2}} 
\def\bold3{{\bf  3}}
\def\bold4{{\bf 4}}
\def\bold5{{\bf 5}}
\def\bold6{{\bf 6}}
\def\bold7{{\bf 7}}
\def\bold8{{\bf 8}}
\def\bold9{{\bf 9}}

\def\P2Skly{\PP^2_{Skly}}

\def\End{\operatorname {End}}

\def\Hom{\operatorname {Hom}}

\def\op{{\operatorname {op}}}

\def\End{\operatorname{End}}

\def\Hom{\operatorname{Hom}}

\def\id{\operatorname{id}}

\def\Mod{{\sf Mod}}

\def\op{\operatorname{op}}

\def\Spec{\operatorname{Spec}}

\def\ul1{\operatorname{\underline{1}}}

\def\G{\mathop{\underline{\underline{\it \Gamma}}}\nolimits}

\def\l{\leftarrow}

\def\a{\alpha}

\def\c{\gamma}

\def\l{\lambda}

\def\G{\Gamma}
\def\L{\Lambda}

\def\sA{{\sf A}}
\def\sB{{\sf B}}
\def\sC{{\sf C}}

\def\cal{\mathcal}

\def\cE{{\cal E}}
\def\cF{{\cal F}}
\def\cG{{\cal G}}

\def\cN{{\cal N}}
\def\cO{{\cal O}}

\def\Qcoh{{\sf Qcoh}}

\def\dirlim{\mathop{\vtop{\baselineskip -100pt\lineskip -1pt\lineskiplimit 0pt
\setbox0\hbox{lim}\copy0\hbox to \wd0{\rightarrowfill}}}\limits}
\def\invlim{\mathop{\vtop{\baselineskip -100pt\lineskip -1pt\lineskiplimit 0pt
\setbox0\hbox{lim}\copy0\hbox to \wd0{\leftarrowfill}}}\limits}

\def\I11{{1 \kern -0.8pt \! \mbox{l}}}
\def\mumu{{\mu\kern-4.2pt\mu}}
\def\bfmu{{\mu\kern-4.2pt\mu}}
\def\2slash{\backslash \! \backslash}

\def\boxtimes{\setbox0\hbox{$\Box$}\copy0\kern-\wd0\hbox{$\times$}}

\begin{document}

\title[A generalization of Watts's Theorem]{A generalization of Watts's Theorem: \\
Right exact functors on module categories}
\author{A. Nyman and S. P. Smith}
\address{Department of Mathematics, University of Montana, Missoula, MT 59812-0864}
\email{NymanA@mso.umt.edu} 
\address{Department of Mathematics, University of
Washington, Box 354350, Seattle, Washington 98195}
\email{smith@math.washington.edu}
\date{\today}
\thanks{The first author was partially supported by the National Security Agency under grant NSA H98230-05-1-0021. The second author was supported by NSF grant DMS-0245724.}

\maketitle

\begin{abstract}
Watts's Theorem says that a right exact functor $F:\Mod R \to \Mod
S$ that commutes with direct sums is isomorphic to $-\otimes_R B$
where $B$ is the $R$-$S$-bimodule $FR$. The main result in this
paper is the following: if $\sA$ is a cocomplete abelian category and
$F:\Mod R \to \sA$ is a right exact functor commuting with direct
sums, then $F$ is isomorphic to $- \otimes_R \cF$ where $\cF$ is a
suitable $R$-module in $\sA$, i.e., a pair $(\cF,\rho)$ consisting
of an object $\cF \in \sA$ and a ring homomorphism $\rho:R \to
\Hom_\sA(\cF,\cF)$.  Part of the point is to give meaning to the
notation $-\otimes_R \cF$. That is done in the paper by Artin and
Zhang \cite{AZ} on Abstract Hilbert Schemes. The present paper is
a natural extension of some of the ideas in the first part of
their paper.
\end{abstract}

\pagenumbering{arabic}

\section{Introduction}
Let $R$ and $S$ be rings and let $\Mod R$ and $\Mod S$ denote the
category of right $R$-modules and right $S$-modules, respectively.
Watts's Theorem, which was proved by Eilenberg \cite{E} and Watts \cite{W} at about the same 
time,  is the following:

\begin{thm} \label{thm.classical}
Suppose $F:\Mod R \to \Mod S$ is a right exact functor commuting
with direct limits.  Then $F \cong -\otimes_R B$ where $B$ is an
$R$-$S$-bimodule.
\end{thm}

Let $\sB(\Mod R,\Mod S)$ denote the full subcategory of the
category of functors from $\Mod R$ to $\Mod S$ consisting of right
exact functors commuting with direct limits.   The next result is a slightly more precise
version of Theorem  \ref{thm.classical}.

\begin{thm} \label{thm.generalization}
The functor $\Psi:\Mod (R^{\op} \otimes_\mathbb{Z} S) \to \sB(\Mod R,\Mod
S)$ induced by the assignment $B \mapsto -\otimes_R B$ is an
equivalence of categories.
\end{thm}
Theorem \ref{thm.classical} is then just the fact that the functor
$\Psi$ is essentially surjective.

The main result of this paper (Theorem \ref{thm.watts}) is that if
$\Mod S$ is replaced by an arbitrary cocomplete\footnote{An additive category is {\sf cocomplete} 
if it has arbitrary direct sums. This is Grothendieck's condition Ab3.}
 category $\sA$,
then a version of Theorem \ref{thm.generalization} still holds.
One of the obvious hurdles in proving such a theorem is to have a
sensible notion of tensor product in this context.  We use the
tensor product functor that was defined in \cite[Thm. 3.7.1]{Po}
and investigated in detail in  \cite{AZ} (see Section
\ref{subsection.tensor}).

In Proposition \ref{prop.scheme}, we specialize our main result to
the case that $\sA$ is the category of quasi-coherent sheaves on a
scheme $Y$. This version of the main result is used extensively in
\cite{NS} to prove a structure theorem for right exact functors
between categories of quasi-coherent sheaves on schemes.

\section{Preliminaries}
Throughout this paper, $k$ is a fixed commutative ring, $R$ is a
$k$-algebra, and $\gamma:k \to R$ is the homomorphism giving $R$
its $k$-algebra structure.

\subsection{$k$-linearity}
 Let  $\sA$ be an additive category. We say $\sA$ is
 {\it $k$-linear} if for all objects $X$ and $Y$ in $\sA$,
$\Hom_\sA(X,Y)$ is a $k$-module and composition of morphisms
is $k$-bilinear.  Equivalently,  $\sA$ is $k$-linear if there is a ring homomorphism
$$
c: k \to \End(\id_\sA)
$$
from $k$ to the ring of natural transformations from the identity functor to itself.

The first definition tells us that for each object $X \in \sA$ and each $a \in k$ there is a morphism
$a_X:X \to X$ such that
\begin{equation} \label{eqn.lincompat}
a_Y \circ f =f \circ a_X
\end{equation}
 for all $a \in k$ and $f \in \Hom_\sA(X,Y)$.
The second definition tells us there are natural transformations
$c(a):\id_\sA \to \id_\sA$ for each $a \in k$, and therefore
associated morphisms $c(a)_X:X \to X$ for each $a \in k$ and $X \in \sA$.
The connection between the two definitions is that $$c(a)_X = a_X$$ for all $a \in k$ and $X$ in
$\sA$.

The $k$-linear structure on $\Mod R$ is given by
\begin{equation}
\label{ModR.k-linear}
a_M(m)= m.\c(a).
\end{equation}
for all $M \in \Mod R$, $m \in M$, and $a \in k$.

\subsection{$k$-linear functors}

Let $\sC$ and $\sA$ be $k$-linear categories. A functor $F:\sC \to
\sA$  is  {\it $k$-linear} if the natural maps $\Hom_\sC(X,Y) \to
\Hom_\sA(FX, FY)$ are $k$-linear for all $X$ and $Y$ in $\sC$.
Equivalently, $F$ is $k$-linear if $F$ is additive and
$$
F(a_Y)=a_{FY}
$$
for all $a \in k$ and $Y \in \Mod R$.

 We write
$$
\sB_k(\sC,\sA)
$$
for the full subcategory of the category of functors $\sC \to \sA$ consisting of $k$-linear right exact
functors that commute with direct limits. We use the letter $\sB$ to remind us of
bimodules.

It is surely well known that an adjoint to a $k$-linear functor is
again $k$-linear but we could not find a proof in the literature so provide one for completeness.

\begin{prop}
\label{prop.adj.linear}
Let $\sC$ and $\sA$ be $k$-linear categories.
Let $G:\sA \to \sC$ be a functor having a left adjoint $F$. If $G$ is $k$-linear so is $F$.
\end{prop}
\begin{pf}
Let $X \in \sC$, and let
$$
\nu:\Hom_\sA(FX,FX) \to \Hom_\sC(X,GFX)
$$ be the adjoint isomorphism. By the functoriality of the adjoint isomorphisms the diagrams
$$
\UseComputerModernTips
\xymatrix{
\Hom_\sA(FX,FX) \ar[rr]^\nu  \ar[d]_{- \circ Ff} &&  \Hom_\sC(X,GFX) \ar[d]^{-\circ f}
\\
\Hom_\sA(FX,FX) \ar[rr]_\nu  &&  \Hom_\sC(X,GFX)
}
$$
and
$$
\UseComputerModernTips
\xymatrix{
\Hom_\sA(FX,FX) \ar[rr]^\nu  \ar[d]_{g \circ -} &&  \Hom_\sC(X,GFX) \ar[d]^{Gg \circ -}
\\
\Hom_\sA(FX,FX) \ar[rr]_\nu  &&  \Hom_\sC(X,GFX)
}
$$
commute for all $X$ in $\sC$, all $f \in \Hom_\sC(X,X)$, and all $g \in \Hom_{\sA}(FX,FX)$.

Let $\theta \in \Hom_\sA(FX,FX)$ be an element in the top left corner of the diagrams.
Let $f=a_X$ and $g=a_{FX}$. The commutativity therefore gives
\begin{align*}
\nu(\theta \circ F(a_X)) & = \nu(\theta) \circ a_X \qquad \hbox{and}
\\
\nu(a_{FX} \circ \theta) & =G(a_{FX}) \circ \nu(\theta).
\end{align*}
But $\nu(\theta):X \to GFX$ is a $k$-linear morphism so $ \nu(\theta) \circ a_X  = a_{GFX} \circ \nu(\theta)$.
Since $G$ is $k$-linear, $G(a_{FX}) = a_{GFX}$. Hence
$$
\nu(\theta \circ F(a_X))  =a_{GFX} \circ \nu(\theta) = G(a_{FX})  \circ \nu(\theta) =  \nu(a_{FX} \circ \theta).
$$
But $\nu$ is an isomorphism so
$$
\theta \circ F(a_X)  = a_{FX} \circ \theta.
$$
Now take $\theta=\id_{FX}$ to get $F(a_X)=a_{FX}$, so showing that $F$ is $k$-linear.
\end{pf}

\subsection{The category $\sA_R$}

For the remainder of this paper, we let $\sA$ denote a $k$-linear
cocomplete category.

A  {\it left $R$-module in $\sA$} is a pair $(\cF,\rho)$ where $\cF$ is an object in $\sA$
and $\rho:R \to \End_\sA \cF$ is a $k$-algebra homomorphism.
Popescu \cite[p. 108]{Po} calls $(\cF,\rho)$ a left $R$-object of $\sA$.
Let $(\cF,\rho)$ and $(\cG,\rho')$ be left $R$-modules in $\sA$. We define the set of {\it $R$-module maps} from  $(\cF,\rho)$ to $(\cG,\rho')$ to be
$$
\Hom_R(\cF,\cG): = \bigl\{ \a \in \Hom_\sA(\cF,\cG) \; | \; \rho'(r) \circ \a=\a \circ \rho(r) \;
\hbox{ for all } \; r \in R\bigr\}.
$$
Using these $R$-module maps as morphisms we then obtain a category $\sA_R$,
 the category of left $R$-modules in $\sA$.

Suppose $(\cF,\rho) \in \sA_R$. If $\cG \in \sA$, then
$\Hom_\sA(\cF,\cG)$ becomes a right $R$-module through the
composition map
$$
\Hom_\sA(\cF,\cG) \times \Hom_\sA(\cF,\cF) \to \Hom_\sA(\cF,\cG),
$$
i.e.,
$$
\a.r:= \a \circ \rho(r)
$$
for $\a \in \Hom_\sA(\cF,\cG)$ and $r \in R$. This allows
us to view $\Hom_\sA(\cF,-)$ as a functor $\sA \to \Mod R$.

\smallskip

For each $x \in R$,  let  $\mu_x:R \to R$ be the right $R$-module
homomorphism  $\mu_x(r):=xr$.

\begin{lemma} \label{lemma.ar}
Suppose $F \in \sB_k(\Mod R,\sA)$. Define the ring homomorphism
$$
\rho:R \to \End_{\sA}FR, \qquad \rho(x):=F(\mu_x).
$$
Then $(FR,\rho) \in \sA_R$.
\end{lemma}

\begin{pf}
To prove the lemma it suffices to show that $\rho$ is a
$k$-algebra homomorphism, i.e., that $(\rho \circ \c)(a) = a_{FR}$
for all $a \in k$. But
$$
\rho(\c(a)) = F(\mu_{\c(a)}) = F(a_R) = a_{FR},
$$
where the second equality is due to (\ref{ModR.k-linear}). Hence the result.
\end{pf}

\subsection{The functor $-\otimes_R \cF$}
\label{subsection.tensor} 
Recall the standing hypothesis that $\sA$ is cocomplete.

Let $(\cF,\rho) \in \sA_R$. By \cite[p.
108]{Po}, the functor $\Hom_\sA(\cF,-): \sA \to \Mod R$ has a left
adjoint.\footnote{It is essential that $\sA$ be cocomplete
for $-\otimes_R \cF$ to exist. For example, if  $R=\ZZ$ and $\sA$ consists of finitely generated abelian 
groups and $\cF=\ZZ$, there is no adjoint. But the hypothesis of cocompleteness is absent from
\cite[p.108]{Po} and parts of  \cite{AZ}.}
  We fix a left adjoint and denote it by $-\otimes_R \cF$.  By
\cite[Proposition B3.1]{AZ}, the functor $-\otimes_R \cF$ is
unique up to isomorphism (of functors) such that
\begin{itemize}
\item{} $R \otimes_R \cF \cong \cF$, and

\item{} $-\otimes_R \cF$ is right exact and commutes with direct
sums.
\end{itemize}

Since the functor $\Hom_\sA(\cF,-)$ is $k$-linear for all $\cF \in
\sA$, Proposition \ref{prop.adj.linear} implies the following:

\begin{cor}
 \label{cor.klinear}
If $(\cF, \rho) \in \sA_R$, then $-\otimes_R \cF$ is $k$-linear.
\end{cor}

\section{The generalization of Watts's Theorem}

\begin{thm} \label{thm.watts}
The functor
$$
\Psi:\sA_R \to \sB_{k}({\sf Mod}R,\sA)
$$
induced by the assignment
$$
\qquad \Psi(\cF) = - \otimes_R \cF,
$$
is an equivalence of categories.
\end{thm}

\subsection{The proof that $\Psi$ is essentially surjective}

\begin{prop}\footnote{After we finished writing this paper we learned that a special case of this result had already
been proved by Brzezinski and Wisbauer \cite[39.3, p.410]{BW} under the hypothesis that the objects 
of $\sA$ are abelian groups.}
\label{prop.general.watts} Let $F \in \sB_{k}({\sf Mod}R,\sA)$.
Then $F \cong -\otimes_R \cF$ where $\cF=FR$.
\end{prop}

\begin{pf} 
Let $\theta_M : M \to \Hom_\sA(\cF,FM)$ be the composition
$$
\UseComputerModernTips
\xymatrix{
M \ar[r]^>>>>>{\L_M } & \Hom_R(R,M) \ar[r]^F & \Hom_{\sA}(\cF,FM)
}
$$
where $\L_M$ is the canonical isomorphism $m \to \l_m$ where
$\l_m(r):=mr$ for all $r \in R$.

Let
$$
\Theta_M:M \otimes_R \cF \to FM
$$
be the map that corresponds to $\theta_M$ under the adjoint isomorphism
$$
\Hom_R(M,\Hom_\sA(\cF,FM)) \cong \Hom_{\sA}(M \otimes_R \cF, FM).
$$

We will show that the $\Theta_M$s define a natural transformation, i.e., if $f:M \to N$
is a homomorphism of right $R$-modules, then the diagram
\begin{equation}
\label{nat.trans?}
\UseComputerModernTips
\xymatrix{
M   \otimes_R \cF\ar[r]^{f \otimes \cF} \ar[d]_{\Theta_M}& N  \otimes_R \cF \ar[d]^{\Theta_N}
\\
FM \ar[r]_{Ff} & FN
}
\end{equation}
commutes.
Define $\eta: \Hom_{\sA}(\cF,FM) \to  \Hom_{\sA}(\cF,FN)$ by $\eta(g):=Ff \circ g$. The left and right  squares in the diagram
$$
\UseComputerModernTips
\xymatrix{
M \ar[r]^>>>>>{\L_M } \ar[d]_f  & \Hom_R(R,M) \ar[r]^F \ar[d] & \Hom_{\sA}(\cF,FM) \ar[d]^\eta
\\
N \ar[r]_>>>>>{\L_N } & \Hom_R(R,N) \ar[r]_F & \Hom_{\sA}(\cF,FN)
}
$$
commute, so $\eta \circ \theta_M=\theta_N \circ f$.

We now consider the diagram
$$
\UseComputerModernTips
\xymatrix{
\Hom(M,\Hom(\cF,FM)) \ar[d] \ar[rr]^\sim &  &  \Hom(M\otimes \cF,FM)  \ar[d]
\\
\Hom(M,\Hom(\cF,FN))  \ar[rr]^\sim &&  \Hom(M \otimes \cF,FN)&&
\\
 \Hom(N,\Hom(\cF,FN))  \ar[u] \ar[rr]^\sim  &&  \Hom(N \otimes \cF,FN),   \ar[u]
}
$$
whose verticals are induced by $f$ and whose horizontals are the
adjoint isomorphism.  The top and bottom rectangles of this
diagram commute by the functoriality of the adjoint isomorphisms.
The maps $\theta_M$ and $\theta_N$ belong to the top and bottom
Hom-sets of the left-hand column and their images in $
\Hom(M,\Hom(\cF,FN))$ are the same because $\eta \circ \theta_M
=\theta_N \circ f$. It follows that the images of $\Theta_M$ and
$\Theta_N$ in $\Hom(M \otimes \cF,FN)$ are the same. In other
words,
$$
Ff \circ \Theta_M = \Theta_N \circ (f \otimes \cF)
$$
which proves that (\ref{nat.trans?}) commutes and hence that the $\Theta_M$s define a natural transformation
$$
\Theta: - \otimes_R \cF \to F.
$$

Because $(F \circ \L_R)(x)=F(\mu_x) = \rho(x)$, $\theta_R:R \to
\Hom_\sA(\cF,\cF)$ is the map giving $\cF$ its $R$-module
structure, so the corresponding map  $\Theta_R : R \otimes_R \cF
\to \cF$ is an isomorphism. Since the functors $-\otimes_R \cF$
and $F$ commute with direct sums, $\Theta_M$ is an isomorphism for
all free $R$-modules $M$. Since  $-\otimes_R \cF$ and $F$ are
right exact it follows that $\Theta_M$ is an isomorphism whenever
$M$ is the cokernel of a map between free $R$-modules. But every
$R$-module is of that form so $\Theta_M$ is an isomorphism for all
$M$. Hence $\Theta$ is an isomorphism of functors.\footnote{The argument in the last part of the proof is a result of B. Mitchell. See \cite[39.1, p.409]{BW} for more details.} 
\end{pf}

Proposition \ref{prop.general.watts} says that the
functor $\Psi$ in Theorem \ref{thm.watts} is essentially
surjective.

\subsection{$R \otimes_R -:\sA_R \to \sA_R$ is isomorphic to the identity functor}

Let $(\cF,\rho) \in \sA_R$ and let $\cN \in \sA$. The composition
\begin{equation}
\label{defn.theta}
\UseComputerModernTips
\xymatrix{
\Hom_{\sA}(R \otimes_R \cF,\cN) \ar[r]^{\sim} &  \Hom_R(R,\Hom_\sA(\cF,\cN)) \ar[r]^\sim
& \Hom_\sA(\cF,\cN),
}
\end{equation}
where the first map is the adjoint isomorphism and the second is
the canonical isomorphism $\psi \mapsto \psi(1)$, induces an
isomorphism of functors $$\Hom_{\sA}(R \otimes_R \cF,-) \to
\Hom_\sA(\cF,-)$$ which, by the Yoneda Lemma, corresponds to a
unique isomorphism
$$
\UseComputerModernTips
\xymatrix{
\xi_\cF: \;   \cF   \ar[r]^\sim  &   R \otimes_R \cF.
}
$$

The next result is a slightly sharper form of \cite[Prop. B3.1(a)]{AZ}.

\begin{prop}
\label{prop.AZ.B3.1} The diagram
\begin{equation}
\label{eq.theta} 
\UseComputerModernTips \xymatrix{
 R \otimes_R \cF   \ar[r]^{R \otimes \phi} &  R \otimes_R \cG
\\
\cF \ar[r]_{\phi}  \ar[u]^{\xi_\cF}  & \cG \ar[u]_{\xi_\cG}
}
\end{equation}
commutes for all $\cF,\cG \in \sA_R$ and all $\phi \in
\Hom_R(\cF,\cG)$. Therefore, the maps $\xi_\cF$ provide an
isomorphism
$$
\xi:\id_{\sA_R} \longrightarrow (R  \otimes_R -)
$$
of functors.
\end{prop}
\begin{pf}
By the Yoneda lemma, the commutivity of (\ref{eq.theta}) is
equivalent to the condition that for all $\cN \in \sA$ the outer
rectangle in the diagram
 \begin{equation}
 \label{eq.theta.nat.trans}
\UseComputerModernTips
\xymatrix{
\Hom_{\sA}(R \otimes_R \cF,\cN) \ar[d]_{\cong} &&  \ar[ll]_{-\circ(R \otimes \phi)}  \Hom_{\sA}(R \otimes_R \cG,\cN) \ar[d]^{\cong}
\\
\Hom_R(R,\Hom_\sA(\cF,\cN)) \ar[d]_{\cong} &&  \ar[ll]_{\G} \Hom_R(R,\Hom_\sA(\cG,\cN)) \ar[d]^{\cong}
\\
\Hom_\sA(\cF,\cN) && \ar[ll]^{-\circ \phi} \Hom_\sA(\cG,\cN)
}
\end{equation}
commutes, where the vertical arrows are the factorizations  in (\ref{defn.theta}) that are used
to define $\xi_\cF$ and $\xi_\cG$, and
$$
\G(\psi)(x) :=\psi(x) \circ \phi
$$
for all $x \in R$ and $\psi \in \Hom_R(R,\Hom_\sA(\cG,\cN))$.

The uppermost square of (\ref{eq.theta.nat.trans}) commutes by functoriality of
the adjoint isomorphism.  Going clockwise around the lower square,
the image in $\Hom_\sA(\cF,\cN)$ of  $\psi \in \Hom_R(R,\Hom_\sA(\cG,\cN))$ is $\psi(1) \circ \phi$.
Going counter-clockwise around the lower square,
the image of  $\psi$ in $\Hom_\sA(\cF,\cN)$ is $\Gamma(\psi)(1) = \psi(1) \circ \phi$.
Hence the lower square commutes.

It follows that the outer rectangle commutes.
\end{pf}

\subsection{The proof that $\Psi$ is full} 
 Let $\cF$ and $\cG$ be objects in $\sA_R$ and let $\phi \in \Hom_R(\cF,\cG)$. 
 Because (\ref{eq.theta}) commutes  $R \otimes \phi$ is non-zero if $\phi$ is non-zero. Hence
$\Psi$ is faithful.

\subsection{The proof that $\Psi$ is faithful}
To complete the proof of Theorem \ref{thm.watts}, it remains to
show that $\Psi$ is full. To that end,  let $$\tau:-\otimes_R \cF
\to -\otimes_R \cG$$ be a natural transformation. We must show
there is a homomorphism $\phi \in \Hom_R(\cF,\cG)$  such that
$\tau_M=M \otimes \phi$ for all $M \in \Mod R$.

Define
$$
\phi:= \xi_\cG^{-1} \circ \tau_R \circ \xi_{\cF}.
$$
It follows from the commutativity of (\ref{eq.theta}) that $R \otimes
\phi=\tau_R$. By Lemma \ref{lem.nat.trans} below, it follows that $M
\otimes \phi=\tau_M$ for all $M \in \Mod R$. In other words,
$\Psi(\phi)=\tau$.

\begin{lem}
\label{lem.nat.trans}
Let  $\sC$ be an abelian category and let $F,G:\Mod R \to \sC$ be right exact functors that commute
with direct sums. Let $\tau,\tau':F \to G$ be natural transformations. If $\tau_R=\tau'_R$, then
$\tau=\tau'$.
\end{lem}
\begin{pf}
Let $M_i$, $i \in I$, be a collection of right $R$-modules. Then there is a natural map
$$
\bigoplus_{i \in I} FM_i \to F \Bigl(\bigoplus_{i \in I} M_i\Bigr)
$$
and the fact that $F$ commutes with direct sums says that this map
is an isomorphism. By the universal property of colimits, there is
a commutative diagram
$$
\UseComputerModernTips
\xymatrix{
\bigoplus_{i \in I} FM_i \ar[r]   \ar[d]_{\oplus \tau_{M_i}} & F \Bigl(\bigoplus_{i \in I} M_i\Bigr) \ar[d]^{\tau_{\oplus M_i}}
\\
\bigoplus_{i \in I} GM_i \ar[r]  & G \Bigl(\bigoplus_{i \in I}
M_i\Bigr). }
$$
Since the horizontal maps are isomorphisms,  if $\tau_{M_i}=\tau'_{M_i}$ for all $i$, then
$$
\tau_{\oplus M_i}=\tau'_{\oplus M_i}.
$$
In particular, it follows that $\tau_P=\tau'_P$ for all free $R$-modules $P$.

Let $M$ be a right $R$-module and let $P \to Q \to M \to 0$ be an exact sequence in which $P$ and
$Q$ are free $R$-modules.  Then there is a commutative diagram
$$
\UseComputerModernTips
\xymatrix{
FP \ar[d]_{\tau_P} \ar[r] & FQ  \ar[d]_{\tau_Q} \ar[r] & FM   \ar[r] & 0
\\
GP   \ar[r] & GQ   \ar[r] & GM    \ar[r] & 0,
}
$$
 and a unique map $FM \to GM$ making the diagram commute, namely $\tau_M$.
 Since $\tau_P=\tau'_P$ and $\tau_Q=\tau'_Q$, it follows that $\tau_M=\tau'_M$.
\end{pf}

\section{An application}
Throughout this section, let $X$ denote a $k$-scheme.  If
$X=\operatorname{Spec }R$, we let
$$
\widetilde{(-)}:{\sf Mod }R \to {\sf Qcoh }X
$$
be the quasi-inverse to the global sections functor defined in
\cite[II, Definition, p. 110]{Hart}.

\begin{example}
Let $f:Y \to X$ be a morphism from an arbitrary scheme to an
affine scheme $X = \Spec R$.  Then $f^* \circ \widetilde{(-)}:\Mod
R \to \Qcoh Y$ is a right exact functor commuting with direct
sums. Proposition \ref{prop.general.watts} says that $f^* \circ
\widetilde{(-)} \cong - \otimes_R \cO_Y$ where $\cO_Y$ is made
into an $R$-module via the ring homomorphism
$$
R \to \Hom_R(R,R) \to
\Hom_Y(f^*\cO_X,f^*\cO_X) \to \Hom_Y(\cO_Y,\cO_Y)
$$
where the first map sends $r \in R$ to multiplication by $r$, the
second map is induced by $f^{*}\circ \widetilde{(-)}$ and the
third isomorphism is induced by the natural isomorphism
$f^{*}\cO_X \cong O_Y$.
\end{example}

The motivation for this present paper lies in our paper \cite{NS}.
There we consider $k$-schemes $X$ and $Y$ and $k$-linear
functors  $F:\Qcoh X \to \Qcoh Y$ that are right exact and commute with direct sums.
One source of such functors is the following. Let $\cF$ be a quasi-coherent sheaf on $X \times_k Y$, and
define
\begin{equation}
-\otimes_{\mathcal{O}_{X}}\cF :=
\operatorname{pr}_{2*}(\operatorname{pr}_{1}^{*}(-)
\otimes_{\mathcal{O}_{X \times_{k} Y}} \cF)
\end{equation}
where $\operatorname{pr}_{i}:X \times_{k} Y \to X,Y$, $i=1,2$, are
the obvious projections.  

Two warnings are necessary. First, a functor of the form
$-\otimes_{\mathcal{O}_{X}}\cF$ is not always in
$\sB_{k}({\sf Qcoh }X,{\sf Qcoh }Y)$.  This happens, for example,
if $Y=\Spec k$, $X=\PP^1_k$ and
$\cF=\mathcal{O}_{X \times Y} \cong
\Gamma(X,-)$.

Second, an object in $\sB_{k}({\sf Qcoh }X,{\sf Qcoh
}Y)$ is not always isomorphic to one of the form
$-\otimes_{\mathcal{O}_{X}}\cF$.  For example, when
$Y=\Spec k$ and $X=\PP^1_k$, the functor $H^1(X,-)$ is not of this form 
\cite[Proposition  5.4]{NS}.

The question motivating \cite{NS} is whether $F$ is isomorphic to
a functor of the form $-\otimes_{\mathcal{O}_{X}}\cF$.  It follows
from Theorem \ref{thm.watts} that this is always the case if $X$
is affine, as we now show.

\begin{prop}
 \label{prop.scheme}
 Let $R$ be a commutative $k$-algebra and  $Y$ a $k$-scheme.
 Write $X := \operatorname{Spec }R$.
Then the inclusion functor
$$
{\sf Qcoh }(X \times_{k} Y ) \to \sB_{k}({\sf Qcoh }X,{\sf Qcoh }Y),
\qquad
\cF \mapsto
-\otimes_{\mathcal{O}_{X}}\cF,
$$
is an equivalence of  categories.
\end{prop}
\begin{pf}
By \cite[II, exercise 5.17e]{Hart}, the functor
$$
\operatorname{pr}_{2*}:{\sf Qcoh }(X \times_{k} Y )\to {\sf Qcoh
}(\operatorname{pr}_{2*}\mathcal{O}_{X \times_{k}Y})
$$
is an equivalence, where ${\sf Qcoh
}(\operatorname{pr}_{2*}\mathcal{O}_{X \times_{k}Y})$ denotes the
category of quasi-coherent $\mathcal{O}_{Y}$-modules with
$\operatorname{pr}_{2*}\mathcal{O}_{X \times_{k}Y}$-module
structure.  Furthermore, it is straightforward to check that the
functor
$$
{\sf Qcoh }(\operatorname{pr}_{2*}\mathcal{O}_{X \times_{k}Y}) \to
({\sf Qcoh }Y)_R
$$
induced by the assignment $\cE \mapsto (\cE,\rho)$, where $\rho:R
\to \Hom_Y(\cE,\cE)$ is defined through the
$\operatorname{pr}_{2*}\mathcal{O}_{X \times_{k}Y}$-structure of
$\cE$, is an equivalence.  By Theorem \ref{thm.watts}, the functor
$$
({\sf Qcoh }Y)_R \to {\sf B}_{k}(\Mod R,{\sf Qcoh }Y)
$$
induced by the assignment $(\cE,\rho) \mapsto -\otimes_R \cE$ is
an equivalence.  Therefore, the functor
$$
{\sf Qcoh }(X \times_{k} Y) \to {\sf B}_{k}(\Mod R,{\sf Qcoh }Y)
$$
induced by the assignment $\cF \mapsto -\otimes_R
\operatorname{pr}_{2*}\cF$ is an equivalence.  By the uniqueness
properties of the functor $-\otimes_R \cE$ described in Section
\ref{subsection.tensor}, we have an isomorphism of functors
$$
\UseComputerModernTips
\xymatrix{
-\otimes_R  \operatorname{pr}_{2*}\cF
\ar[rr]^{\sim}    &&
\widetilde{(-)} \otimes_{\mathcal{O}_{X}}\cF
}
$$
in ${\sf B}_{k}(R,{\sf Qcoh }Y)$.  It follows that the functor
$$
{\sf Qcoh }(X \times_{k} Y) \to {\sf B}_{k}(\Mod R,{\sf Qcoh }Y)
$$
induced by the assignment $\cF \mapsto \widetilde{(-)}
\otimes_{\mathcal{O}_{X}}\cF$ is an equivalence.  The claim
follows easily from this.
\end{pf}


\begin{thebibliography}{11}


\bibitem{AZ}
M. Artin and J.J. Zhang,
Abstract Hilbert Schemes,
{\it Alg. and Repn. Theory,} {\bf 4} (2001) 305-394.

\bibitem{BW}
T. Brzezinski and R. Wisbauer, {\it Corings and Comodules,}
Lond. Math. Soc. Lect. Note Ser. 309, Camb. Univ. Press, 2003.

\bibitem{E}
S. Eilenberg, Abstract description of some basic functors,
{\it J. INd. Math. Soc.,} {\bf 24} (1960) 231-234.

\bibitem{Hart}
R. Hartshorne, {\it Algebraic Geometry}, Springer-Verlag, New York
1977.

\bibitem{NS}
A. Nyman and S.P. Smith, Watts's Theorem for Schemes, {\it in
preparation}.

\bibitem{Po}
N. Popescu, {\it Abelian Categories with Applications to Rings and
Modules}, Academic Press, London 1973.

\bibitem{W}
C. E. Watts, Intrinsic characterizations of some additive functors,
{\it Proc. Amer. Math. Soc.,} {\bf 11} (1960) 5-8.

\end{thebibliography}
\end{document}